\documentclass[10pt,a4paper,twoside]{amsart}

\usepackage{amsfonts,amssymb,euscript,eufrak}
\usepackage{amscd,amsmath}
\usepackage{cite}
\usepackage[all]{xy}

\newcommand{\Q}{\mathbb{Q}}

\newcommand{\Z}{\mathbb{Z}}

\newcommand{\al}{\alpha}

\newcommand{\ol}{o}

\newcommand{\G}{\mathbb{G}}

\newcommand{\hato}{\,\hat{\otimes}}

\newcommand{\car}{\stackrel{\cong}{\longrightarrow}}

\newcommand{\C}{\hat{C}}
\newcommand{\dual}{k[\epsilon]}
\newcommand{\PM}{\mathfrak{PM}}
\newcommand{\M}{\mathfrak{M}}\newcommand{\m}{\mathfrak{m}}
\newcommand{\limi}{\underleftarrow{\lim}}
\newcommand{\ModA}{{\rm Mod}_G^{\rm pro\,aug}(A)^{\rm fl}}
\newcommand{\Modk}{{\rm Mod}_G^{\rm pro\,aug}(k)^{\rm fl}}
\newcommand{\ModB}{{\rm Mod}_G^{\rm pro\,aug}(B)^{\rm fl}}

\newcommand{\ModAo}{{\rm Mod}_G^{\rm pro\,aug}(A)}
\newcommand{\Modko}{{\rm Mod}_G^{\rm pro\,aug}(k)}
\newcommand{\ModBo}{{\rm Mod}_G^{\rm pro\,aug}(B)}

\newcommand{\Modsmk}{{\rm Mod}_G^{\rm sm}(k)}

\newcommand{\Irr}{{\rm Mod}_{G,0}^{\rm pro\,aug}(k)}

\newcommand{\Nr}{N}

\newcommand{\EndA}{{\rm End}_A^{cont}}

\newcommand{\Endk}{{\rm End}_k^{cont}}

\newcommand{\Def}{\mathcal{D}}

\newcommand{\oP}{\overline{P}}
\newcommand{\ochi}{\overline{\chi}}
\newcommand{\oep}{\overline{\epsilon}}
\newcommand{\oV}{\bar{V}}

\newcommand{\af}{\mathfrak{a}}
\newcommand{\bfr}{\mathfrak{b}}
\newcommand{\mf}{\mathfrak{m}}

\newcommand{\Bangu}{{\rm Ban}_G(K)^{\leq 1}}
\newcommand{\BanG}{{\rm Ban}_G(K)}

\newenvironment{pr}{\it Proof:\rm}{\hfill $\Box$\newline}
\newtheorem{theo}{Theorem}[section]
\newtheorem{prop}[theo]{Proposition}
\newtheorem{cor}[theo]{Corollary}

\newtheorem{lem}[theo]{Lemma}

\begin{document}
\title[Deformations]{On unitary deformations of smooth modular representations}

\author[Tobias Schmidt]{Tobias Schmidt}
\address{Mathematisches Institut\\ Westf\"alische Wilhelms-Universit\"at
M\"unster\\ Einsteinstr. 62\\ D-48149 M\"unster, Germany}
\email{toschmid@math.uni-muenster.de}

\maketitle
\begin{abstract}
Let $G$ be a locally $\Q_p$-analytic group and $K$ a finite
extension of $\Q_p$ with residue field $k$. Adapting a strategy of
B. Mazur (cf. \cite{MazurI}) we use deformation theory to study
the possible liftings of a given smooth $G$-representation $\rho$
over $k$ to unitary $G$-Banach space representations over $K$. The
main result proves the existence of a universal deformation space
in case $\rho$ admits only scalar endomorphisms. As an application
we let $G={\rm GL}_2(\mathbb{Q}_p)$ and compute the fibers of the
reduction map in principal series representations.
\end{abstract}
\section{Introduction}

Let $G$ a locally $\Q_p$-analytic group and $K$ a finite extension
of $\Q_p$ with ring of integers $\ol$ and residue field $k$. The
aim of the present note is to study the set of possible liftings
of a given smooth $G$-representation $\rho$ over $k$ to unitary
$G$-Banach space representations over $K$. To do this we adapt the
techniques of deformation theory for representations of profinite
groups as developed by B. Mazur (cf. \cite{MazurI}) to our present
situation.  We prove, in case $\rho$ admits only scalar
endomorphisms (e.g. admissible and absolutely irreducible), the
existence of a formal scheme ${\rm Spf}\; R(G,K,\rho)$ over $\ol$
which depends functorially on the datum $(G,K,\rho)$ and respects
elementary operations on $\rho$ such as tensor product. Moreover,
its $\ol$-rational points biject canonically with the isomorphism
classes of unitary liftings of $\rho$. The ring $R(G,K,\rho)$ is a
local profinite $\ol$-algebra with residue field $k$ which is
noetherian if and only if the $k$-vector space of extensions ${\rm
Ext}_G^1(\rho,\rho)$ is of finite dimension. Basic features of the
geometry of $Spec \;R(G,K,\rho)$ such as dimensions or the number
of irreducible components remain unclear at this point.

The ring $R(G,K,\rho)$ represents a deformation problem for
Iwasawa modules which is based on the simple observation (due to
V. Paskunas, cf. \cite{Paskunas})) that the duality functor
introduced by P. Schneider and J. Teitelbaum (cf. \cite{STBanach})
on the category of unitary representations is compatible with
Pontryagin duality over the profinite ring $\ol$. By work of M.
Emerton (cf. \cite{EmertonI}) the dual categories admit
generalizations to complete local noetherian $\ol$-algebras which
provide a natural framework to study deformations of (the
Pontryagin dual of) $\rho$.

We emphasize two cases in which our main result is well-known. If
$G$ is a compact group and $\rho$ is admissible absolutely
irreducible our result follows directly from work of B. Mazur. On
the other hand, in the important case of $G={\rm GL}_2(\Q_p)$ our
result was essentially established by M. Kisin building on work of
P. Colmez (cf. \cite{Colmez10}) and M. Schlessinger (cf.
\cite{Schlessinger}). In \cite{KisinDEF} this result is proved to
show the essential surjectivity of Colmez's Montreal functor, an
important result in the $p$-adic Langlands programme for ${\rm
GL}_2(\Q_p)$. In this light we hope our result will prove
important in extending the $p$-adic Langlands programme to other
groups than ${\rm GL}_2(\Q_p)$.

\bigskip

Let us briefly outline the paper. We begin by recalling and
establishing some basic results on pseudocompact rings (sect.
$2$). For sake of clarity and to have a greater flexibility in
future applications we then proceed axiomatically and work in a
quite general setting as follows. Let $k$ be a field (finite or
not) of characteristic $p>0$ and $\rho$ be a smooth
$G$-representation over $k$. Let $\ol$ be an arbitrary complete
local noetherian ring with residue field $k$. We introduce a
category $\C$ of coefficient rings consisting of commutative local
pseudocompact $\ol$-algebras $A$ such that the structure map
$\ol\rightarrow A$ gives an isomorphism on residue fields.
Analogous to work of M. Emerton (loc.cit.) we introduce a suitable
category of Iwasawa modules over such an $A$ and study its basic
properties (sect. $3$). By results of A. Brumer (cf.
\cite{Brumer}) this category depends "naturally" on $A$ which
results in a deformation functor $D_\rho$ for $\rho$. It is
straightforward (sect. $3$) to prove the representability of
$D_\rho$ for characters $\rho$ and to explain the notion of
deformation conditions in our setting. Section $4$ contains the
proof of the representability of $D_\rho$ in case $\rho$ has only
scalar endomorphisms and $k$ is finite and establishes the ${\rm
Ext}^1$-criterion. We remark straightaway that Schlessinger theory
(loc.cit.) is not applicable in our situation since we refrain
from any finiteness assumption on the tangent space of $D_\rho$.
Instead, we proceed directly from A. Grothendieck's fundamental
representability theorem (cf. \cite{GrothendieckDesII}) using
ideas of M. Dickinson (cf. \cite{GouDef}, Appendix 1).
In section 5 we show the usual functorial properties of the
universal deformation ring. This is almost a formality.


In the final section we turn to Banach space representations and
specialize the deformation theory to the situation where $\ol$ is
given as ring of integers in $K$. This yields the relation between
unitary lifts of $\rho$ and rational points of ${\rm
Spf}\;R(G,K,\rho)$. We illustrate this method with an application
to the group $G={\rm GL}_2(\Q_p)$ and compute the set of unitary
deformations in case $\rho$ equals a principal series
representation.  This heavily builds upon results of M. Emerton
concerning his functor of 'ordinary parts' (\cite{EmertonII}). We
also remark that the structure of such principal series
representations $\rho$ over $k$ has been made explicit by the work
of Barthel/Livn\'e (\cite{BarthelLivne}).




~\\{\it Acknowledgements.} The author would like to thank the
Arbeitsgruppe of Peter Schneider at M\"unster University whose
Oberseminar in the winter term 2005/06 provided a motivating
opportunity for the author to learn about deformation theory and
its applications to group representations. He would also like to
thank Matthew Emerton for an interesting conversation on this
work.

\section{Pseudocompact rings}\label{pseudocompactrings}

For any unital ring $A$ we let $\M(A)$ be the abelian category of
left unital $A$-modules. If $A$ is left noetherian then the
finitely generated left $A$-modules form a full subcategory
$\M_{fg}(A)$ of $\M(A)$. A left and right noetherian ring will be
called noetherian.

\bigskip

A {\it left pseudocompact ring} is a complete Hausdorff
topological unital ring $A$ which admits a system of open
neighbourhoods of zero consisting of left ideals $\af$ such that
$A/\af$ has finite length as a left $A$-module (cf. \cite{SGA3}).
In particular, $A$ equals the topological inverse limit of the
artinian quotients $A/\af$ each endowed with the discrete
topology. A morphism of left pseudocompact rings is by definition
a continuous unital ring homomorphism. A left artinian ring with
the discrete topology is evidently left pseudocompact. More
generally, the topology on a left pseudocompact ring $A$ which is
left noetherian is uniquely determined and coincides with the adic
topology defined by the Jacobson radical of $A$.

Let $A$ be a left pseudocompact ring. A complete Hausdorff
topological left unital $A$-module $M$ is called {\it left
pseudocompact} if it has a system of open neighbourhoods of zero
consisting of submodules $M'$ such that $M/M'$ has finite length.
A morphism between two left pseudocompact modules is by definition
a continuous $A$-linear map. It necessarily has closed image.
Borrowing notation from \cite{SV1} we denote the category of left
pseudocompact $A$-modules by $\PM(A)$. It is abelian with exact
projective limits and the forgetful functor
$\PM(A)\rightarrow\M(A)$ is faithful and exact and commutes with
projective limits (cf. \cite{GabrielPhD}, IV.3. Thm. 3,
\cite{GvB}, Prop. 3.3). An arbitrary direct product of left
pseudocompact modules is left pseudocompact in the product
topology. A left pseudocompact module $M$ is called {\it
topologically free} if it is topologically isomorphic to a product
$\prod_I A$ with some index set $I$. The set of images in $M$ of
the "unit vectors" $(...,0,1,0,...)\in\prod_I A$ under such an
isomorphism is called a {\it left pseudobasis} of $M$.
If $A$ is left noetherian then any finitely generated abstract
left $A$-module has a unique left pseudocompact topology. We thus
have a natural fully faithful and exact embedding
$\M_{fg}(A)\rightarrow\PM(A)$.

There are obvious "right" versions of the statements above.

\bigskip

Remark: For sake of clarity we point out the following. In
\cite{Brumer}, A. Brumer defines a pseudocompact ring to be a
complete Hausdorff topological unital ring $A$ that admits a
system of open neighborhoods of zero consisting of two sided
ideals $I$ for which $A/I$ is an Artin ring (i.e. satisfies the
descending chain condition for chains of two-sided ideals). If $A$
is commutative then this is equivalent to $A$ being left and right
pseudocompact.

\bigskip

For the rest of this section let us fix a left pseudocompact ring
$A$ which is {\it commutative}. Evidently, it is then right
pseudocompact and we will not distinguish between left and right
$A$-modules.

Given a pseudocompact module $M$ over such an $A$ write $M^*:={\rm
Hom}(M,A)$ for the $A$-module of morphisms $M\rightarrow A$ in
$\PM(A)$. We obtain a functor
\begin{equation}\label{projective}M\mapsto M^*\end{equation} between $\PM(A)$ and $\M(A)$ that
changes direct products into direct sums. If $A$ is artinian this
establishes an anti-equivalence of categories between projective
objects in $\PM(A)$ and $\M(A)$ respectively (cf. \cite{SGA3},
0.2.2). In general, if $A^{I}:=\prod_I A$ denotes a topologically
free module on a pseudobasis indexed by $I$ and
\[M_I(A):={\rm End}_{\PM(A)}(A^{I})\] its endomorphism ring
we obtain an isomorphism of abstract $A$-modules
\begin{equation}\label{iso}M_I(A)\car\prod_I
(A^{I})^*=\prod_I\oplus_I A,\end{equation} natural in $A$. In the
light of (\ref{iso}) we sometimes view elements of $M_I(A)$ as
infinite '$I\times I$-matrices' with entries from $A$. We always
equip $M_I(A)$ with the compact-open topology. If $A^{I}$ is
locally compact then $M_I(A)$ is a topological ring (\cite{B-GT},
X.\S3.4 Prop.9).

\vskip8pt

If $M,N$ denote two pseudocompact $A$-modules define the
$A$-module \begin{equation}\label{tensor}M\hato_A
N:=\limi_{M',N'}\;M/M'\otimes_A N/N'\end{equation} where $M'$ and
$N'$ run through the open submodules of $M$ and $N$ respectively.
If each $M/M'\otimes_A N/N'$ is endowed with the discrete topology
the projective limit topology makes $M\hato_A N$ a pseudocompact
$A$-module. Indeed, given $M',N'$ as above there exists an open
ideal $\af\subseteq A$ such that $\af M\subseteq M'$ and $\af
N\subseteq N'$
so that $M/M'\otimes_A N/N'$ is a finitely generated module over
the artinian ring $A/\af$ and therefore of finite $A$-length. The
binary operation $\hato_A$ on $\PM(A)$ is associative and
commutative with $A$ as a unit object and functorial in both
variables. It commutes with projective limits and direct products
(cf. \cite{SGA3}, 0.3.5/6). Now let $\phi: A\rightarrow B$ be a
morphism between commutative pseudocompact rings and $M\in\PM(A)$.
Define the $B$-module $M\hato_A \;B$ by the exact analogue of
formula (\ref{tensor}). Arguing similarly as above shows $M\hato_A
B$ to be a pseudocompact $B$-module. We obtain a "base change"
functor
\[\phi^*: \PM(A)\longrightarrow \PM(B)\]
which commutes with tensor products, projective limits and direct
products (cf. \cite{SGA3}, 0.5).

\bigskip

A {\it pseudocompact algebra over} $A$ is a topological unital
$A$-algebra $B$ (commutative or not) which admits a system of open
neighbourhoods of zero consisting of twosided ideals $\bfr$ such
that $B/\bfr$ has finite length as an $A$-module (cf.
\cite{Brumer}). For example, an $A$-algebra which is of finite
length as $A$-module is evidently a pseudocompact $A$-algebra in
the discrete topology. The following simple fact will prove useful
in the sequel.

\begin{lem}\label{proGab}
Let $B$ be a topological unital $A$-algebra. Then $B$ is a
pseudocompact $A$-algebra if and only if the underlying
topological $A$-module of $B$ is pseudocompact.
\end{lem}
\begin{pr}
If suffices to show the "if" part. Let $P\subseteq B$ be a left
ideal such that $B/P$ is of finite $A$-length and consider the
multiplication map followed by the natural projection
$B\rightarrow B/P$
\[ \varphi: B\times B\longrightarrow B\longrightarrow B/P.\]
Arguing similarly as in \cite{SGA3}, Lem. 0.3.1 we find a left
open submodule $M$ and a right open submodule $N$ of $B$ of finite
colength such that $\varphi(B,N)=\varphi(M,B)=0$. In other words,
$P$ contains the two-sided ideal of $B$ generated by $M\cap N$
which is open.
\end{pr}

A {\it morphism} between two pseudocompact algebras over $A$ is by
definition a continuous unital $A$-algebra homomorphism. As a
consequence of the above lemma the category of pseudocompact
$A$-algebras has projective limits. Moreover, if $B\rightarrow C$
and $B\rightarrow D$ are morphisms between pseudocompact
$A$-algebras then the completed tensor product $C\hato_B\; D$ is a
pseudocompact $A$-algebra in the obvious way. Finally, if $B$ is a
pseudocompact $A$-algebra and $C$ is a pseudocompact $B$-algebra
then evidently $C$ is a pseudocompact $A$-algebra.

\bigskip

\bigskip

We finally point out the following simple construction. Evidently,
a pseudocompact $A$-algebra $B$ is a left and right pseudocompact
ring.
Furthermore, a discrete (topological) $B$-module has finite length
as a $B$-module if and only if it has finite length as an
$A$-module.
Letting $\phi: A\rightarrow B$ denote the structure map we
therefore have a natural faithful and exact forgetful functor
\begin{equation}\label{scalarres}
\phi_*: \PM(B)\longrightarrow\PM(A).\end{equation}

\section{The deformation functor}

\subsection{Completed group algebras}

Let $A$ be a commutative pseudocompact ring and $H$ be a profinite
group. Writing $\mathcal{N}$ for the system of open normal
subgroups of $H$ we denote by
\[A[[H]]:=\limi_{N\in\mathcal{N}} \;A[H/N]\]
the completed group algebra of $H$ over $A$. It is a pseudocompact
$A$-algebra with respect to the projective limit topology and the
correspondance \[H\mapsto A[[H]]\] is a covariant functor from
profinite groups to pseudocompact $A$-algebras (cf. \cite{Brumer},
Sect. 2). The anti-involution $h\mapsto h^{-1}$ on $H$ identifies
$A[[H]]$ and $A[[H]]^{opp}$ as pseudocompact $A$-algebras making
it unnecessary to distinguish between left and right modules. Let
$\PM(A[[H]])$ be the abelian category of (left) pseudocompact
$A[[H]]$-modules.
\bigskip

Remark:  Suppose $A$ is noetherian and $H$ is locally
$\mathbb{Q}_p$-analytic. Since the pseudocompact topology on $A$
is defined by the Jacobson radical (cf. sect. 2) a mild
generalization of [Emeb], Thm. 2.1.1 shows that $A[[H]]$ is
noetherian. We shall make no use of this fact.

\bigskip

Recall that a topological $A$-module $M$ is called a (topological)
{\it left $H$-module} over $A$ if it has a $H$-action by
$A$-linear maps such that the map
\[ H\times M\longrightarrow M\]
giving the action is continuous. It follows that a discrete
$A$-module of finite length is an $H$-module if and only if it is
a (discrete) $A[[H]]$-module. A projective limit argument shows
that a pseudocompact $A$-module $M$ is an $H$-module if and only
if it is a (pseudocompact) $A[[H]]$-module (cf. \cite{Brumer}, p.
454/455).

After these preliminaries let $M,N\in\PM(A[[H]])$ be given. The
diagonal $H$-action on the pseudocompact $A$-module $M\hato_A N$
is then continuous and, by our initial remarks, extends therefore
to a pseudocompact $A[[H]]$-module structure.
The resulting binary operation $\hato_A$ on $\PM(A[[H]])$ is
associative, commutative and functorial in both variables. The
usual augmentation homomorphism $A[[H]]\rightarrow A$ provides a
unit object.

Now consider a pseudocompact $A$-algebra $B$ and let $\phi:
A\rightarrow B$ be the structure map. The base change $\phi^*:
\PM(A)\rightarrow\PM(B)$ commutes with projective limits (cf.
sect. $2$). Thus, the compatible system of natural isomorphisms
\[B\otimes_A A[H/N]\car B[H/N],~~~N\in\mathcal{N}\] induces a
natural isomorphism $\phi^*(A[[H]])\car B[[H]]$ which is
multiplicative. This discussion yields a functor
\begin{equation}\label{nummer}\phi^*_H: \PM(A[[H]])\rightarrow
\PM(B[[H]])\end{equation} compatible with $\phi^*$ via the
appropriate forgetful functors. Given another pseudocompact
$B$-algebra $B'$ with structure map $\psi$ we evidently have
\begin{equation}\label{assH}(\psi\circ\phi)^*_H=\psi^*_H\circ\phi^*_H.\end{equation}
Let additionally, $\PM(A[[H]])^{\rm fl}$ denote the full
subcategory of $\PM(A[[H]])$ consisting of modules on underlying
topologically free $A$-modules (similarly for $B$). Since the
tensor product commutes with direct products (cf. sect. $2$) the
functor $\phi^*_H$ is seen to respect these subcategories.

\subsection{Augmented representations}
Let $\ol$ be a fixed commutative complete local noetherian ring
and $A$ a commutative pseudocompact $\ol$-algebra. We now bring in
a locally $\Q_p$-analytic group $G$ and introduce a certain
category of $G$-representations over $A$.

\bigskip
Let $A[G]$ be the group algebra of $G$ over $A$. A {\it
$G$-representation over $A$} is simply a (left) $A[G]$-module.
Following \cite{EmertonI} such a representation $M$ is called {\it
augmented} if the induced $A[H]$-action extends to an
$A[[H]]$-action on $M$ for every
compact open subgroup $H$ of $G$. 
Let $M$ be an augmented $G$-representation on a pseudocompact
$A$-module. It is called a \textit{pseudocompact augmented}
$G$-representation over $A$ if the topology on $M$ makes it a
pseudocompact $A[[H]]$-module for every compact open subgroup $H$
of $G$. By the above this is equivalent to requiring that the
induced $H$-action makes $M$ into an $H$-module for every compact
open subgroup $H$ of $G$.

\begin{lem}\label{contG}
On a pseudocompact augmented $G$-representation $M$ over $A$ the
group $G$ acts by continuous automorphisms.
\end{lem}
\begin{pr}
Pick a compact open subgroup $H$ of $G$. Invoking the duality
$M\mapsto M^\vee$ between $\PM(A[[H]])$ and the category of
discrete topological $A[[H]]$-modules (cf. \cite{Brumer}) it
suffices to consider $M^\vee$ with the $G$-action induced by
functoriality. But $M^\vee$ is a discrete $H$-module and so the
$G$-action is even jointly continuous on $M^\vee$.
\end{pr}

We define a morphism between two such representations to be a
$G$-equivariant morphism in $\PM(A)$. Since $A[H]\subseteq A[[H]]$
is dense it is plain that any such morphism extends to a morphism
in $\PM(A[[H]])$ for every compact open subgroup $H$ of $G$.
Borrowing notation from \cite{EmertonI} we denote the resulting
category by $\ModAo.$
\begin{lem}\label{abelian}
The category $\ModAo$ is an $A$-linear abelian tensor category.
\end{lem}
\begin{pr}
Let $H$ be a compact open subgroup of $G$. Since $\PM(A)$ and
$\PM(A[[H]])$ are abelian categories $\ModAo$ is abelian as well.
Furthermore, given $M,N\in\PM(A[[H]])$ the diagonal $G$-action on
\[M\hato_A N\in\PM(A[[H]])\] for any compact open subgroup $H\subseteq G$ makes the latter module
a pseudocompact augmented $G$-representation over $A$. This yields
the desired tensor product on $\ModAo$.
\end{pr}

Remark: In case of a compact discrete valuation ring $\ol$ and a
complete local noetherian ring $A$ with finite residue field the
category $\ModAo$ was introduced by M. Emerton. It is Pontryagin
dual to certain smooth $G$-representations over $A$ and plays a
central role in M. Emerton's theory of ordinary parts of
admissible representations. For more details we refer to
\cite{EmertonI} and \cite{EmertonII}.

\subsection{Functors on coefficient algebras}
Let as before $\ol$ be a commutative complete local noetherian
ring. We now define two subcategories of commutative pseudocompact
$\ol$-algebras which will serve as coefficient algebras within the
upcoming deformation theory.

Let $\mf$ be the maximal ideal of $\ol$ with residue field $k$.
Let $\hat{C}$ be the full subcategory of pseudocompact algebras A
over $\ol$ that are commutative local rings and such that the
structure map $\ol\rightarrow A$ is local and induces an
isomorphism on residue fields. Let $C$ denote the full subcategory
of $\C$ consisting of discrete algebras having finite length as
$\ol$-module. Without recalling the precise definition of a {\it
pro-object} (cf. \cite{GrothendieckDesII}, A.2) we have the
following simple observation.

\begin{lem}\label{pro} The category of pro-objects of
$C$ is equivalent to $\C$.
\end{lem}
\begin{pr}
Let us denote by $Pro(.)$ the passage from a category to the its
category of pro-objects and pro-morphisms. Let $C'$ be the
category of all discrete $\ol$-algebras having finite length as
$\ol$-module. Mapping a pseudocompact $\ol$-algebra to the system
of all its artinian quotients induces an equivalence between
pseudocompact $\ol$-algebras and $Pro(C')$ (cf.
\cite{GrothendieckDesII}, A.5). Restricting this functor to $\C$
yields a fully faithful functor into $Pro(C)$. Given an element
$(R_i)_i$ in $Pro(C)$ the projective limit $\limi_i R_i$ lies in
$\C$ and hence, this functor is essentially surjective.
\end{pr}

We now bring in a set-valued covariant functor on $C$
\[D:C\longrightarrow Sets.\] The category $C$ contains $k$ as a
terminal object and admits finite products and finite fiber
products (cf. \cite{MazurII}, Lem. IV.\S14). As to the latter,
recall that if $\phi_i: A_i\rightarrow A_0$ are two morphisms in
$C$ their fiber product is given as the equalizer
\[ A_1\times_{A_0}A_2=\{(a_1,a_2)\in A_1\times A_2: \phi_1(a_1)=\phi_2(a_2)\}\]
with ring structure induced from $A_1\times A_2$. In this
situation $D$ is called {\it left exact} if it respects finite
products and finite fiber products. Furthermore, since $\C$
identifies with the pro-objects of $C$, the functor $D$ being {\it
pro-representable} is tantamount to being of the form ${\rm
Hom}_{\C}(R,.)$ with some $R\in\C$ (cf. \cite{GrothendieckDesII},
A.2).

Let $k[\epsilon]=k[x]/x^2$ be the ring of dual numbers viewed as
an object in $C$. If $D$ is pro-representable the set
$D(k[\epsilon])$ evidently has a natural $k$-vector space
structure (the "tangent space" of $D$).
\begin{theo}\label{grothendieck}
The functor $D: C\rightarrow Sets$ is pro-representable if and
only if it is left exact. In this situation the representing ring
$R$ is noetherian if and only if the $k$-vector space
$D(k[\epsilon])$ has finite dimension $d$. In this case, $R$
equals a quotient of the formal power series ring
$\ol[[x_1,...,x_d]]$.
\end{theo}
\begin{pr}
This follows directly from A. Grothendieck's fundamental
representability theorem (cf. \cite{GrothendieckDesII}, Prop.
A.3.1/A.5.1).
\end{pr}

Suppose we now have a functor $D: \C\rightarrow Sets$ on the
larger category $\C$. By the above lemma $\C$ is stable under
arbitrary projective limits.
The above discussion therefore shows that $D$ is representable as
a functor on $\C$ if and only if it commutes with projective
limits and the restriction of $D$ to $C$ is pro-representable.

\subsection{Deformations}
We define the deformation problem and state the main
representability result. We keep the assumptions of the previous
subsection but {\bf assume} additionally that the residue field of
$\ol$ has characteristic $p>0$.

\bigskip

Given a pseudocompact $A$-algebra $B$ let $\phi: A\rightarrow B$
be the structure map. For every $H\in\mathcal{N}$ we have the base
change $\phi^*_H$ commuting with forgetful functors (cf.
(\ref{nummer})) and, hence, a functor
\begin{equation}\label{basechange}\phi^*_G: \ModAo\longrightarrow\ModBo\end{equation} respecting the
full subcategories $\ModA$ and $\ModB$ of modules which are
topologically free over $A$ and $B$ respectively. Given another
pseudocompact $B$-algebra $B'$ with structure map $\psi$ one has
\begin{equation}\label{assG}(\psi\circ\phi)^*_G=\psi^*_G\circ\phi^*_G\end{equation} according to
(\ref{assH}). After these preliminaries we fix once and for all an
element
\[N\in\Modko.\] 
Let $I$ be an index set of a pseudobasis for the topologically
free $k$-module underlying $N$. Invoking the duality $N\mapsto
N^*$ (cf. (\ref{projective})) we see that the cardinality $|I|$
does not depend on the choice of pseudobasis. Given a local
pseudocompact $\ol$-algebra $A\in\C$ with residue homomorphism
$\phi: A\rightarrow k$ we consider couples $(M,\al)$ such that
$M\in\ModA$ and
\[\al: k\hato_A M=\phi^*_G(M)\car N\] is an isomorphism in
$\Modk$.
\begin{lem}
Given $M\in \PM(A)$ the natural map
\[k\hato_A M\longrightarrow M/\overline{\m M}\]
is an isomorphism in $\PM(k)$. If $\m$ is finitely generated, the
subset $\m M\subseteq M$ is closed.
\end{lem}
\begin{pr}
The first statement is \cite{SGA3}, 0.3.2. Suppose $\m$ is
finitely generated. The $A$-submodule $\m M\subseteq M$ equals the
image of a morphism in $\PM(A)$ of the form $\prod M\rightarrow M$
where the product is indexed by finitely many generators of the
ideal $\m$. It is therefore closed.
\end{pr}

A morphism of couples \[(M,\al)\longrightarrow (M',\al')\] is by
definition a morphism $M\rightarrow M'$ in $\ModA$ such that the
resulting diagram
\begin{equation}\label{mor}\xymatrix{
M/\overline{\m M} \ar[r]^-{\sim}_-{\al}\ar[d]& N \ar[d]^{=}\\
  M'/\overline{\m M'}  \ar[r]^-{\sim}_-{\al'}& N}
\end{equation}
is commutative. We denote the set of isomorphism classes of such
couples by $D_\Nr(A)$. As usual elements in $D_\Nr(A)$ will be
called {\it deformations} of $N$ to $A$ and we will often
abbreviate $M$ for $[M,\al]$ when no confusion can arise. Since
base change to $k$ commutes with arbitrary direct products any
pseudobasis of $M\in D_\Nr(A)$ must have cardinality $I$, too.
This shows
\begin{lem}\label{under}
If $[M,\al]\in D_\Nr(A)$ then $M\simeq A^{I}$ in $\PM(A)$.
\end{lem}
By associativity
(\ref{assG}) of the base change $\phi^*_G$ we obtain a covariant
set-valued functor
\[D_\Nr: \C\longrightarrow Sets\] such that $D_N(k)$ is a
singleton. By Lem. \ref{abelian} the category $\Modko$ is abelian
whence we have the $k$-vector space ${\rm Ext}^1(N,N)$ of {\it
Yoneda extensions} of $N$ by itself.
\begin{prop}\label{tangent2}
Let $k$ be finite. There is a natural bijection
\[D_\Nr(k[\epsilon])\car {\rm Ext}^1(N,N)\]
which is $k$-linear in case $D_\Nr$ is representable.
\end{prop}
\begin{pr}
This is a standard phenomenon in deformation theory (cf.
\cite{MazurII}, Prop. V.\S22) and most probably true without the
restriction on the field $k$. In any case, let $A=\dual$ with
$\pi: A\rightarrow k$ and $\iota: k\rightarrow A$ the canonical
maps. Let $[M,\al]\in D_\Nr(A)$. Invoking an isomorphism $M\simeq
A^{I}$ in $\PM(A)$ (cf. Lem. \ref{under}) we have $\epsilon
M\simeq\prod_Ik\simeq N$ in $\PM(k)$. Since multiplication by
$\epsilon$ is a homeomorphism onto its image the restriction map
induced by $\epsilon M\subseteq M$ and $\iota$
\[res: \EndA(M)\longrightarrow\Endk(\epsilon M)\]
is continuous for compact-open topologies. 
The diagram
\[\xymatrix{
  A[[H]] \ar[r]^-{cont} & {\rm End}^{cont}_A(M)\ar[r]^-{res} & \Endk(\epsilon M)\\
  A[H] \ar[r]^{\subseteq} \ar[u]^{\subseteq} &  A[G]\ar[u] \ar[ur]&.
  }
\]
is commutative and factores in an obvious sense through the
projection $\pi$ applied to coefficients. By the discussion at the
beginning of section \ref{mainproof} below we therefore see that
$\epsilon M\in\Modko$ and that the isomorphism $\epsilon M\simeq
N$ in $\PM(k)$ lifts to an isomorphism in $\Modko$. In other words
\[0\rightarrow \epsilon M\rightarrow M\rightarrow M/\epsilon
M\rightarrow 0\] yields an element in ${\rm Ext}^1(N,N)$. We
obtain a map $D_\Nr(A)\rightarrow{\rm Ext}^1(N,N)$. To construct
the inverse let
\begin{equation}\label{sequence} 0\longrightarrow
N\stackrel{\iota}
{\longrightarrow}M\stackrel{\pi}{\longrightarrow} N\longrightarrow
0\end{equation} be an extension in $\Modko$. By topological
freeness it splits in $\PM(k)$ and we may assume that $M$ admits a
neighbourhood basis of zero consisting of $k$-vector spaces
preserved under the endomorphism $\iota\circ\pi$.
Setting $\epsilon.m:=(\iota\circ\pi)(m)$ makes $M$ a linearly
topologized $A$-module whence $M\in\PM(A)$. An easy argument shows
that the $A$-module on the dual $k$-vector space $M^*$ is free.
Applying the quasi-inverse to (\ref{projective}) we find $M$ to be
topologically free over $A$. Finally, since $\iota, \pi$ are
linear with respect to $A[[H]]$ and $A[G]$ it follows that
$M\in\ModAo$.
This construction gives the inverse. The assertion about linearity
follows immediately from the construction.
\end{pr}

We come to the main result of this section which will be proved in
section \ref{mainproof}.
\begin{theo}\label{main}Let $k$ be finite.
If ${\rm End}_{\Modko}(N)=k$ then $D_{\Nr}$ is representable.
\end{theo}

Remarks: 1. Suppose $N$ is finitely generated over $k[[H]]$ for
some (equivalently any) compact open $H\subseteq G$. Then any
$M\in D_{\Nr}(A)$ is finitely generated over $A[[H]]$. Indeed,
lifting finitely many $k[[H]]$-module generators of
$M/\overline{\m M}=N$ to $M$ we obtain a map $\prod
A[[H]]\rightarrow M$ whose cokernel $Q$ satisfies $Q/\overline{\m
Q}=0$. Since $A$ is local $Q=0$ by the Nakayama lemma (cf.
\cite{SGA3}, 0.3.3).\vskip8pt

2. The theorem is most probably true without the restriction on
$k$. We impose this restriction since our proof makes crucial use
of the compact-open topology. This topology is only well-behaved
for locally compact spaces (\cite{B-GT}, X.\S3.4).\vskip8pt

3. Let $k$ be finite. A pseudocompact $\ol$-algebra $A$ is then a
profinite $\ol$-algebra and a pseudocompact $A$-module is a
profinite $A$-module. In the light of the main result of this work
we could therefore have worked from the beginning on with
profinite algebras and profinite modules. However, to have more
flexibility in future applications it seemed advantageous to us to
produce as much as possible of this theory in the more general
'pseudocompact' language.

\begin{cor} \label{tangent} Let $k$ be finite (so that $D_N$ is representable).
The representing $\ol$-algebra $R$ is noetherian if and only if
$d:= {\rm dim}_k\, {\rm Ext}^1(N,N)<\infty$. In this situation $R$
equals a quotient of the formal power series ring
$\ol[[x_1,...,x_d]]$.
\end{cor}
\begin{pr}
This follows from prop. \ref{tangent2} and thm.
\ref{grothendieck}.
\end{pr}

We give an important example in which the theorem applies. Assume
$\ol$ equals the integers in a finite extension $K/\Q_p$. In
particular, $k$ is a {\it finite} field of characteristic $p>0$.
As usual, a \textit{smooth} $G$-representation $(V,\rho)$ over $k$
is a $k$-vector space $V$ with a $G$-action such that the
stabilizer of each vector $v\in V$ is open in $G$. With
$G$-equivariant linear maps such representations form an abelian
category $\Modsmk$. A full subcategory is formed by
\textit{admissible} $G$-representations $V$ having the property
that the $k$-vector space of $H$-fixed vectors $V^H$ is finite
dimensional for every compact open subgroup $H\subseteq G$.
Regarding a smooth $G$-representation as a discrete torsion
$\ol$-module with a continuous $G$-action Pontryagin duality
\[V\mapsto V^\vee:={\rm Hom}^{cont}_\ol(V,K/\ol)\] over the
profinite ring $\ol$ induces an anti-equivalence of $k$-linear
categories
\begin{equation}\label{catI} (.)^\vee:
\Modsmk\car \Modko.\end{equation} (cf. \cite{EmertonI}, Lem.
2.2.6). It is compatible with base extension relative to finite
field extensions $k\subseteq k'$.

\begin{lem}\label{schur1}
Let $(V,\rho)\in\Modsmk$ be admissible and absolutely irreducible.
Up to a finite extension of $k$ we then have ${\rm End}_{\Modsmk}(\rho)=k$ and therefore $N=\rho^\vee$
satisfies the conditions of the theorem.
\end{lem}
\begin{pr}
This is Schur's lemma for smooth mod $p$ representations (e.g.
\cite{BreuilNYC}, Lem. 4.1 in the case $G=GL_2(\Q_p)$). Let $f$ be
an endomophism of $\rho$. For a fixed open pro-$p$ subgroup $H$ of
$G$ the space $V^H$ is finite dimensional, nonzero (cf.
\cite{Wilson}, Lem. 11.1.1) and stabilized by $f$. After possibly
extending scalars any nonzero eigenvector of $f$ in $V^H$
generates the $G$-representation $V$.
\end{pr}

Remark: Keeping the assumptions on $\ol$ let $A\in\hat{C}$ be
noetherian. In \cite{EmertonI}, 2.2 M. Emerton defines a new
notion of \textit{smooth} $G$-representation over $A$ (which
coincides with the usual one if $A$ is artinian). With morphisms
being $G$-equivariant $A$-linear maps such representations form an
abelian category which is Pontryagin dual to $\ModAo$. In this
situation, the deformation theory above translates therefore
completely to such smooth representations. We will not use this
point of view in this work and therefore refrain from formulating
a precise picture.

\subsection{The compact case}

Suppose $G$ is compact and $k$ is a {\it finite} field of
characteristic $p>0$. We indicate how our deformation theory
reduces to the situation studied by B. Mazur in the seminal
article \cite{MazurI}.

\bigskip

Given an absolutely irreducible object $N\in\Modko$ let $I$ be an
indexing set for a pseudobasis of $N$. Since any open normal
pro-$p$ subgroup of $G$ acts trivially on the smooth
representation $\rho=N^\vee$ (\cite{Wilson}, Lem. 11.1.1) we have
$n:=|I|<\infty$ and $\rho$ is evidently admissible. Furthermore,
any deformation of $N$ to $A\in\C$ has a finite free underlying
$A$-module (cf. Lem. \ref{under}) and thus, $D_{\Nr}$ describes
the equivalence classes of continuous lifts
\[G\longrightarrow {\rm GL}_n(A)\]
of $N$ to $A$.

Let $M_n(k)$ be the group of $n\times n$-matrices with entries
from $k$. Letting $g.M:=\rho(g)M\rho(g)^{-1}$ for $M\in M_n(k),
g\in G$ makes $M_n(k)$ a (discrete) $G$-module which we denote by
${\rm Ad}(\rho)$ (the so-called \textit{adjoint representation} of
$\rho$). The group ${\rm GL}_n(k[\epsilon])$ is the semidirect
product of the groups $1+\epsilon M_n(k)$ and ${\rm GL}_n(k)$. Let
$\rho_1$ be a deformation of $\rho$ to $k[\epsilon]$. We have a
map
$$c(\rho_1):  G\rightarrow {\rm Ad}(\rho): g\mapsto m_g$$ where $(1+\epsilon
m_g)\rho(g):=\rho_1(g)$. It is a $1$-cocycle inducing an
isomorphism of $k$-vector spaces
$$D_N(k[\epsilon])\car H^1(G,{\rm Ad}(\rho))$$ (cf. \cite{MazurI},
p.399). The \textit{finiteness condition} (in the sense of
[loc.cit.], 1.1) is satisfied since $G$ contains an open pro-$p$
subgroup of finite rank (cf. \cite{DDMS}, Cor. 8.33). We therefore
have $d_1:={\rm dim}_k H^1(G,{\rm Ad}(\rho))<\infty$ (cf.
\cite{MazurII}, Prop. \S21.2a). In this situtation our theorem
(\ref{main}) is therefore equivalent to the existence statement
\cite{MazurI}, Prop.1. In particular, the deformation ring
$R=R(\rho)$ is a (noetherian) quotient of
$\ol[[x_1,...,x_{d_1}]]$. Let $d_2:={\rm dim}_k H^2(G,{\rm
Ad}(\rho))$ and $\m\subset\ol$ the maximal ideal. We have

$${\rm Krull~dim} (R/\m R)\geq d_1-d_2$$

with equality in case $d_2=0$ (the \textit{unobstructed case}); in
this latter case $R=\ol[[x_1,...,x_{d_1}]]$ ([loc.cit.], Prop.2).

\vskip8pt

Remark: It is conjectured that in fact ${\rm Krull~ dim} (R/\m
R)=d_1-d_2$ (compare \cite{GouDef}, p. 287). In \cite{Boeckle}
this conjecture is proved in many cases. In case $G$ equals a
global Galois group it has been pointed out by B. Mazur that this
conjecture should be viewed as a generalization of
\textit{Leopoldt's conjecture} (cf. \cite{MazurI}, 1.10).

\subsection{The case of a character}\label{character}
Let $k$ be finite of characteristic $p>0$ and $N$ such that
$|I|=1$. Under a reasonable technical assumption we will compute
the universal deformation ring $R$ and the universal deformation.
As in the compact case $R$ will not depend on the particular
choice of such an $N$ (cf. Prop. \ref{twisting} below).

Recall that any commutative local pseudocompact ring is {\it
henselian} (cf. \cite{Nagata}, Thm. V.30.3). Suppose $A$ is such a
ring with maximal ideal $\m$ and finite residue field $k=A/\m$. We
have a short exact sequence
\[1\longrightarrow 1+\m\longrightarrow A^\times\longrightarrow
k^\times\longrightarrow 1\] which is canonically split. Let $s:
k^\times\rightarrow\ol^\times$ be the splitting in the case
$A=\ol$.

We assume in the following that the Hausdorff abelianization
$\bar{G}:=G/\overline{[G,G]}$ of $G$ is topologically finitely
generated. Since there are only finitely many continuous group
homomorphisms $G\rightarrow k[\epsilon]^\times$ the universal
deformation ring is noetherian. By \cite{EmertonA}, Prop. 6.4.1
the inclusion of the maximal compact open subgroup $\bar{G}_0$
into $\bar{G}$ induces a (noncanonical) isomorphism
$\bar{G}_0\times\Z^r\cong \bar{G}$ of locally $\Q_p$-analytic
groups for some unique integer $r\geq 0$. If $\Gamma$ denotes the
pro-$p$ completion functor (cf. \cite{Ribes}, 3.2) we have
$\Gamma(\bar{G})=\Z^r_p\times \Gamma(\bar{G}_0)$. The canonical
homomorphism $\bar{G}\rightarrow\Gamma(\bar{G})$ is therefore
continuous with respect to the quotient topology on $\bar{G}$ and
hence, so is the composed homomorphism
\begin{equation}\label{composed} G\rightarrow \bar{G}\rightarrow
\Gamma(\bar{G}).\end{equation} Let $R:=\ol[[\Gamma(\bar{G})]]$. It
is a profinite local noetherian $\ol$-algebra in $\C$ endowed with
a continuous homomorphism $\Gamma(\bar{G})\rightarrow R^\times$.

After these preliminaries let $\bar{\chi}: G\rightarrow k^\times$
be the continuous homomorphism describing the $G$-action on $N$.
Composing $\bar{\chi}$ and the map (\ref{composed}) with the
splitting $s$ and the homomorphism $\Gamma(\bar{G})\rightarrow
R^\times$ respectively we obtain two continuous homomorphisms
$\chi_0: G\rightarrow \ol^\times$ and $\gamma: G\rightarrow
R^\times$.
\begin{prop}
The homomorphism
\[\chi^{univ}=\chi_0\cdot \gamma: G\longrightarrow R^\times\]
equals the universal deformation.
\end{prop}
\begin{pr}
The following argument is a generalization of \cite{GouDef}, Prop.
3.13. Fix a compact open subgroup $H\subseteq G$. First note that
$H\stackrel{\bar{\chi}}{\longrightarrow} k^\times\rightarrow
\ol^\times$ is a continuous homomorphism between profinite groups
whence a continuous algebra homomorphism $R[[H]]\rightarrow
R[[\ol^\times]]$. The homomorphism $\ol^\times\rightarrow
R^\times$ induced by the algebra structure of $R$ is continuous
whence a continuous algebra homomorphism
$R[[\ol^\times]]\rightarrow R$ by the universal property of
\ol[[.]] applied to the profinite ring $R$.
The composite
\[R[[H]]\longrightarrow R[[\ol^\times]]\longrightarrow R\]
then coincides with $\chi_0$ on the subring $R[H]$. Secondly, the
continuous map $\gamma$ restricted to $H$ has image in
$1+\m_R\subseteq R^\times=k^\times\times (1+\m_R)$ since
$\Gamma(\bar{G})$ is a pro-$p$ group. On the one hand, this yields
a continuous algebra homomorphism
\[R[[H]]\longrightarrow R[[1+m_R]]\longrightarrow R\]
that coincides with $\gamma$ on the subring $R[H]$. We therefore
see $\chi^{univ}\in{\rm Mod}_G^{\rm pro\,aug}(R)^{\rm fl}$. On the
other hand, it shows that $\chi^{univ}$ is a deformation of
$\bar{\chi}$.

Now suppose $[\chi]$ is a deformation of $\bar{\chi}$ to a
noetherian ring $A\in \C.$ Since $1+m_A$ is an abelian pro-$p$
group the continuous homomorphism $\chi_0^{-1}\cdot\chi:
G\rightarrow 1+m_A$ factores through
$G\rightarrow\Gamma(\bar{G})$. Indeed, it evidently factores
through $G\rightarrow\bar{G}$ inducing two continuous maps into
$1+m_A$: on the source $\Z^r$ with the discrete topology and on
the source $\bar{G}_0$ with the induced topology as open subgroup
of $\bar{G}$. By a straightforward argument these maps remain
continuous when giving both sources the topology coming from the
families of subgroups of finite $p$-power index. This yields the
claim.

All in all we obtain a continuous map
\[f_\chi: R=\ol[[\Gamma(\bar{G})]]\longrightarrow
\ol[[1+m_A]]\longrightarrow A\] that specializes $\chi^{univ}$ to
$\chi$. It follows that $\chi^{univ}$ is the universal
deformation.
\end{pr}

We give an example in which our hypothesis on $\bar{G}$ is
satisfied.
\begin{prop}\label{noetherian} Let $\G$ denote a connected reductive group
over $\Q_p$ and let $G$ denote the group of its $\Q_p$-rational
points. Then $\bar{G}$ is topologically finitely
generated.
\end{prop}
\begin{pr}
If $Z$ denotes the center of $G$ the natural homomorphism $Z\times
[G,G]\rightarrow G$ has finite kernel and cokernel. 
Being a torus $Z$ is
topologically finitely generated and, hence, so is $\bar{G}$.
\end{pr}


Example: If $G={\rm GL}_n(\Q_p)$ we have
\[\bar{G}=G/[G,G]=\Q_p^\times=p^{\Z}\times U\times\mu \] where
$U=1+p^\kappa\Z_p$ with $\kappa=2$ if $p=2$ and $1$ otherwise and
where $\mu$ denotes the subgroup of roots of unity in
$\Q_p^\times$. The group $U$ is topologically isomorphic to $\Z_p$
whence
\[\Gamma(\bar{G})=\Z_p\times\Z_p\times\mu'\] with $\mu'$ the group
with $\kappa$ elements. Thus $R=\ol[[x_1,x_2]]\otimes_\ol
\ol[\mu']$.

\subsection{Deformation conditions}
It is a formality that the usual formalism of {\it deformation
conditions} works in our setting. We suppose in the following that
the deformation functor $D_\Nr$ is representable.

Given $A\in\C$ assume that some elements of $D_\Nr(A)$ have been
designated to be "of type $\mathcal{P}$" and that this property is
preserved under the base change $D_\Nr(\phi)$ associated to
morphisms $\phi: A\rightarrow B$ in $\C$. We obtain a subfunctor
\[\mathcal{D}_\Nr\subseteq D_\Nr\] by putting $\Def_\Nr(A):=\{M\in D_\Nr(A):
M {\rm~of~type~}\mathcal{P}\}$ for $A\in\C$.

\begin{prop}
The following conditions are necessary and sufficient for the
representability of $\Def_\Nr$:
\begin{itemize}

\item[(1)] $N\in D_N(k)$ has property $\mathcal{P}$.

\item[(2)] Given a diagram $A_1\rightarrow A_0\leftarrow A_2$ in
$C$, any deformation of $N$ to the fiber product
$A_1\times_{A_0}A_2$ whose base changes to $A_1$ and $A_2$ are of
type $\mathcal{P}$ is of type $\mathcal{P}$.

\item[(3)] If $A\in\C$ is an inverse limit of objects $A_i$ in $C$
and the basechange to $A_i$ of a deformation $M$ of $N$ to $A$ is
of type $\mathcal{P}$ for each $i$ then $M$ is of type
$\mathcal{P}$.
\end{itemize}
\end{prop}
\begin{pr}
Granting the representability of $D_N$ the first two conditions
are tantamount to the fact that $\Def_\Nr$ is left-exact and the
third condition asserts that $\Def_\Nr$ preserves arbitrary
inverse limits. The claim is therefore a direct consequence of
Thm. \ref{grothendieck}.
\end{pr}

Let us give two examples of prominent deformation conditions.
Write $Z\subseteq G$ for the center of $G$. Let $A\in C$. We say a
deformation $M$ of $N$ to $A$ {\it has a central character} if $Z$
acts on $M$ via a group homomorphism $Z\rightarrow A^\times$. If
$Z=\Q_p^{\times}$ so that $p\in Z$ we say a deformation $M$ of $N$
to $A$ has {\it uniformizer acting trivially} if $p$ acts
trivially on $M$. We denote these properties by $\mathcal{P}_i,
i=1,2$ respectively.
\begin{prop}
The deformation functor corresponding to $\mathcal{P}_i$ is
representable if and only if $N$ is of type $\mathcal{P}_i$.
\end{prop}
\begin{pr}
We treat the case $\mathcal{P}_1$ of a central character. The
other case is similar. By definition of the $G$-action
$\mathcal{P}_1$ is preserved under base change and it remains to
show that property (1) implies (2)-(3) above. Let
$A_3=A_1\times_{A_0}A_2, M\in D_N(A_3)$ and $\chi_i$ the central
character of the base change $M_i\in D_N(A_i)$. By topological
freeness and compatibility between base change and direct product
we are reduced to $M=A_3$.
Recalling that $A_3$ equals the equalizer in $A_1\times A_2$ of
the maps $A_i\rightarrow A_0$ the character $\chi_1\times\chi_2:
Z\rightarrow A_1^\times\times A_2^\times$ is seen to take values
in $A_3^\times$. Hence, we obtain (2) and (3) follows by a similar
argument.
\end{pr}

\section{Proof of the main result}\label{mainproof}
We fix an object $A\in \C$ and a compact open subgroup $H\subseteq
G$. We identify once and for all $N\simeq k^{I}$ in $\PM(k)$ by
means of a pseudobasis for $N$. As already observed the
cardinality $|I|$ of such a basis is an invariant of $N$. Let $k$
be finite.

Recall the topological ring \[M_I(A):={\rm End}_{\PM(A)}(A^{I})\]
endowed with the compact-open topology. According to \cite{B-GT},
X.\S3.4 Thm. 3 a jointly continuous action
\[A[[H]]\times A^{I}\longrightarrow A^{I}\]
is the same as a continuous $A$-algebra homomorphism
$A[[H]]\rightarrow M_I(A)$.
Having this in mind we follow a strategy of B. Mazur (cf.
\cite{MazurI}) to rewrite the functor $D_\Nr$ in a more accessible
way. Namely, let $E_A$ denote the set of pseudocompact augmented
$G$-representations on the $A$-module $A^{I}$ that lift $N$. It is
evidently functorial in $A$.

Let ${\rm GL}_I(A):=M_I(A)^\times$ be the group of units in
$M_I(A)$. By what we have just said we may think of an element of
$E_A$ as a commutative diagram
\begin{equation}\label{diag2}\xymatrix{
    A[[H]] \ar[r]^-{cont}  & M_I(A)\\
  A[H] \ar[r]^{\subseteq} \ar[u]^{\subseteq} &  A[G] \ar[u].
  }
\end{equation}
that reduces via $A\rightarrow k$ to the corresponding diagram for
$N$. Note here, that the $G$-action on any $M\in\ModAo$ is
necessarily by continuous automorphisms (cf. Lem. \ref{contG})
whence the right vertical arrow. The group ${\rm GL}_I(A)$ acts on
$M_I(A)$ from the left via conjugation. By acting on the
right-upper corner of diagrams (\ref{diag2}) this induces an
action of the subgroup
\[G_A:=\ker \big{(}{\rm GL}_I(A)\rightarrow {\rm
GL}_I(k)\big{)}=1+\prod_I\oplus_I m_{A}\] on the set $E_A$.
\begin{lem}\label{E_A}
There is a bijection
\[E_A/G_A\car
D_\Nr(A)\] natural in $A$.
\end{lem}
\begin{pr}
Let $\pi: A\rightarrow k$ be the residue homomorphism. The map
$A^{I}\mapsto [A^{I}, 1\otimes \pi^{I}]$ induces an injective map
from $E_A/G_A$ to $D_\Nr(A)$. It is surjective by Lem.
\ref{under}.
\end{pr}

Let us now assume that $N$ has only scalar endomorphisms. Let
$\bar{\rho}$ denote the corresponding element in $E_k$.
Given $\rho\in E_A$ write
\[C(\rho)\subseteq M_I(A)\]
for the $A$-algebra equal to the centralizer in $M_I(A)$ of the
image of $\rho: A[G]\rightarrow M_I(A)$. Recall that a surjection
$\phi: A\rightarrow B$ in $C$ is called a {\it small extension} if
$\ker\phi$ equals a nonzero principal ideal which is annihilated
by the maximal ideal of $A$. It is well-known that every
surjection in $C$ factors into a finite composite of small
extensions (e.g. \cite{GouDef}, Problem 3.1).
\begin{lem}\label{schur}
Let $A\in C$. One has $C(\rho)=A$ for all $\rho\in E_A$.
\end{lem}
\begin{pr}
If $A\rightarrow B$ is a small extension in $C$ with kernel
generated by $t\in A$ and $\rho\in E_A$ then the isomorphism
(\ref{iso}) shows that $\ker \big{(}M_I(A)\rightarrow
M_I(k)\big{)}$ is killed by $t$. Furthermore, $C(\bar{\rho})=k$
since $N$ has only scalar endomorphisms. Taken these facts
together the claim is a straightforward generalization of the
arguments given in the proof of \cite{GouDef}, Lem. 3.8.
\end{pr}

As a corollary the functor $D_\Nr$ is {\it continuous} in the
usual sense:
\begin{cor}\label{cont}
Given $A=\limi_n A_n\in\C$ the natural map
\begin{equation}\label{contmap} D_\Nr(A)\car\limi_n~
D_\Nr(A_n)\end{equation} is a bijection.
\end{cor}
\begin{pr}
By Lem. \ref{pro} we may assume that $A_n$ is an artinian quotient
of $A$ so that $A_{n+1}\rightarrow A_n$ is surjective for all $n$.
It follows that the maps
\begin{equation}\label{sur1}G_{A_{n+1}}\longrightarrow G_{A_n}\end{equation} are surjective for all $n$. To check
surjectivity of the map (\ref{contmap}) let $([M_n,\al_n])_n$ be
an element of the projective
limit. A straightforward argument, using the surjectivity of 
(\ref{sur1}) shows the existence of isomorphisms $\beta_n:
M_{n+1}\otimes_{A_{n+1}}A_n\car M_n$ compatible with the $\al_n$.
Passing to the projective limit using $\limi_n\,A_n[[H]]=A[[H]]$
(and similarly for $A[H], A[G]$) yields a pseudocompact augmented
$G$-representation on $M:=\limi_{\beta_n}\,M_n$. By the lemma
below $M$ is topologically free over $A$ and therefore the desired
preimage. For the injectivity let $M,M'$ be representatives of two
classes in $D_\Nr(A)$ together with isomorphisms $M_n\simeq M'_n$
for all $n$ which are compatible with reductions. Let $\rho,\rho'$
be the corresponding elements in $E_A$. A straightforward argument
using the surjectivity of
\begin{equation}\label{sur2}C(\rho_{n+1})\longrightarrow C(\rho_n)\end{equation}
(Lem. \ref{schur}) shows that we may assume the isomorphisms
$M_n\simeq M'_n$ to be compatible with $A_{n+1}\rightarrow A_n$.
Passage to projective limits yields an isomorphism $M\simeq M'$.
\end{pr}

The following statement was used in the preceding proof.
\begin{lem}\label{flatness}
Let $A=\limi_n A_n\in\C$ with $A_n\in C$ an artinian pseudocompact
quotient of $A$ for all $n$. Let $(M_n)_n$ be a projective system
where each $M_n$ is a pseudocompact topologically free
$A_n$-module. The transition map $M_{n+1}\rightarrow M_n$ is
supposed to be continuous and compatible with $A_{n+1}\rightarrow
A_n$. Then $M:=\limi_n M_n$ equipped with the projective limit
topology is a pseudocompact topologically free $A$-module.
\end{lem}
\begin{pr}
Via the quotient map $A\rightarrow A_n$ we may view each $M_n$ as
a pseudocompact $A$-module. It follows that $M$ is a pseudocompact
$A$-module and,
according to the proof of \cite{SGA3}, Prop. 0.3.7, that $M$ is
topologically free over $A$ if the autofunctor $(.)\hato_A M$ on
$\PM(A)$ is exact. So suppose that
\[\mathcal{E}: 0\longrightarrow N'\longrightarrow
M'\stackrel{\phi}{\longrightarrow} P'\longrightarrow 0\] is a
short exact sequence in $\PM(A)$. Let $\mathfrak{m}_n$ be the
kernel of $A\rightarrow A_n$ and put
$M'_n=\overline{\mathfrak{m}_nM'}$. Since $\mathfrak{m}_n$ is
closed it follows easily from \cite{SGA3}, 0.3.2 that $M'=\limi_n
M'/M'_n$. Putting $N'_n=N\cap M'_n$ and $P'_n=\phi(M'_n)$ yields
the exact sequence
\[\mathcal{E}_n: 0\longrightarrow N'/N_n'\longrightarrow
M'/M'_n\stackrel{\phi}{\longrightarrow} P'/P'_n\longrightarrow 0\]
of artinian $A$-modules for all $n$. Since $\cap M'_n=0$ we have
$\cap N'_n=\cap P'_n=0$ and $\limi_n\mathcal{E}_n=\mathcal{E}$ by
exactness of $\limi_n$ on $\PM(A)$. Since $\hato_A$ commutes with
projective limits we obtain isomorphisms of topological
$A$-modules
\[\mathcal{E}\hato_A M\car \limi_n\; \mathcal{E}_n\hato_A
M_n\car\limi_n\; \mathcal{E}_n\hato_{A_n}M_n.\] Since $M_n$ is
topologically free over $A_n$ for all $n$ the functor $(.)\hato_A
M$ is seen to be exact.
\end{pr}

According to Thm. \ref{grothendieck} and Cor. \ref{cont} $D_\Nr$
is representable if its restriction to $C$ is left exact. Since
$D_\Nr(k)$ is a singleton this reduces to verify that $D_\Nr$
respects fiber products. Let therefore \[A_3=A_1\times_{A_0}A_2\]
be a fiber product in $C$. Writing $E_i:=E_{A_i}$ and
$G_i:=G_{A_i}$ and invoking Lem. \ref{E_A} we have to show that
the natural map of sets
\begin{equation}\label{b}
b: E_3/G_3\rightarrow E_1/G_1\times_{E_0/G_0}E_2/G_2
\end{equation}
is a bijection.

\begin{lem}\label{functortoprings}
The correspondance $A\mapsto M_I(A)$ is a functor from $C$ to
topological rings.

\end{lem}
\begin{pr}
Let us show that $M_I(\phi)$ is continuous at zero for a morphism
$\phi: A\rightarrow B$ in $C$. For an open $B$-submodule $V$ of
$B^{I}$ let $W(V)\subseteq M_I(B)$ be the subset consisting of
homomorphisms with image contained in $V$ (and similar notation
for $B$ replaced by $A$). Since $B^{I}$ is compact it suffices to
see that $M_I(\phi)^{-1}(W(V))$ is open. If $U$ denotes the
preimage in $A^{I}$ of $V$ under $\phi^{I}$ then $W(U)\subseteq
M_I(\phi)^{-1}(W(V))$. This shows that $M_I(\phi)$ is continuous
and the rest is clear.
\end{pr}
\begin{lem}
The natural map
\begin{equation}\label{iso2}M_I(A_3)\car
M_I(A_1)\times_{M_I(A_0)}M_I(A_2)\end{equation} is an isomorphism
of topological rings when the target is equipped with the topology
induced by the direct product topology on $M_I(A_1)\times
M_I(A_2)$.
\end{lem}
\begin{pr}
The map in question, say $\phi$, is certainly a continuous ring
homomorphism. We have a chain of canonical isomorphisms of
abstract groups
\[M_I(A_1)\times_{M_I(A_0)}M_I(A_2)\simeq \prod_I (\oplus_I
A_1)\times_{(\oplus_I A_0)}(\oplus_I A_2)\simeq \prod_I\oplus_I
A_3\simeq M_I(A_3)\] whence $\phi$ is bijective. Since all spaces
$A^{I}_i$ are compact it follows easily from the definition of the
compact-open topology that the inverse $\phi^{-1}$ is continuous.
\end{pr}
\begin{lem}\label{onto}
The map $b$ is surjective.
\end{lem}
\begin{pr}
We adapt an argument of M. Dickinson (cf. \cite{GouDef},Appendix
1) to our situation. For this it will be convenient to think of an
augmented $G$-representation $\rho$ in $E_A, A\in C$ as taking
values $\rho(g)$ in infinite $I\times I$-matrices. We shall
therefore write suggestively $c\rho c^{-1}:=c.\rho$ for $c\in
G_A$. Let \[([\rho],[\sigma])\in E_1/G_1\times_{E_0/G_0}E_2/G_2.\]
Let $m_0$ be the maximal ideal of $A_0$. Since $A_0$ is artinian
we have $m_0^n=0$ for some $n\geq 1$. We first prove by induction
on $n$ that $[\rho]$ and $[\sigma]$ have representatives in $E_1$
and $E_2$ respectively whose images coincide in $E_0$.

To do this let $\phi: A_1\rightarrow A_0$ and $\psi:
A_2\rightarrow A_0$ be the transition maps in the fiber product.
By the induction hypothesis we may assume $n>1$ and that
$\phi_*(\rho)=\psi_*(\sigma)$ mod $m_0^{n-1}$. Here, we abbreviate
$\phi_*=E_A(\phi)$ and similarly for $\psi$. Pick an element $c\in
G_0$ such that $c\phi_*(\rho)c^{-1}=\psi_*(\sigma)$ in $E_0$. We
show that there are elements $g\in G_1, h\in G_2$ such that
$\phi_*(g\rho g^{-1})=\psi_*(h\sigma h^{-1})$. This proves the
claim.

Reducing $c$ mod $m_0^{n-1}$ centralizes the image of the
reduction $\phi_*(\rho)$ mod $m_0^{n-1}$ and therefore (Lem.
\ref{schur}) we may assume that $c=1+l$ with $l\in\prod_I\oplus_I
m_0^{n-1}$.
Since $l$ has entries in the finite dimensional $k$-vector space
$m_0^{n-1}$ we may apply {\it mutatis mutandis} \cite{GouDef},
App. 1,Lem. 9.3 and arrive at \[l=\lambda 1+\psi(m_2)-\phi(m_1)\]
with a scalar $\lambda\in m_0^{n-1}$ and $m_i\in \prod_I\oplus_I
A_i$. From now on the claim follows formally from the computations
given in the proof of [loc.cit.], App. 1, Lemma 9.5.

By what we have just shown we may now suppose that
$(\rho,\sigma)\in E_1\times_{E_0} E_2$. The couple gives therefore
rise to a continuous homomorphism
\[ A_1[[H]]\times_{A_0[[H]]}A_2[[H]]\longrightarrow
M_I(A_1)\times_{M_I(A_0)}M_I(A_2)\] compatible with the
$G$-action. Composing it with the obvious continuous ring
homomorphism

$$A_3[[H]]\longrightarrow A_1[[H]]\times_{A_0[[H]]} A_2[[H]]$$

as well as the inverse of the map (\ref{iso2}) yields a preimage
in $E_3$.
\end{pr}

Let again $A_3=A_1\times_{A_0}A_2$ in $C$. For $\rho_i\in E_i$
write
\[G(\rho_i)\subseteq G_i\]
for the stabilizer of $\rho_i$ in $G_i$. According to Lem.
\ref{schur}
\[G(\rho_i)=\big{(}1+\prod_I\oplus_I m_{A_i}\big{)}\cap
C(\rho_i)=1+m_{A_i}.\]
\begin{lem}
The map $b$ is injective.
\end{lem}
\begin{pr}
Assume $b([\rho])=b([\sigma])$ with $[\rho],[\sigma] \in E_3/G_3$.
For $\rho\in E_3$ write $\rho_i$ for the image in $E_i$ and
similarly for $\sigma\in E_3$. Pick $(g_1,g_2)\in G_1\times G_2$
with $\rho_i=g_i.\sigma_i$ in $E_i$. The "top left entry" of the
$I\times I$-matrix of $g_i\in G_i=1+\prod_I\oplus_I m_{A_i}$ lies
in $1+m_{A_i}\subseteq A_i^\times$. Multiplying by a scalar we may
therefore assume this entry is equal to one. Let $\bar{g}_i$
denote the image of $g_i$ in $G_0$. Since $\rho_0=\sigma_0$ we
have $\bar{g}_2^{-1}\bar{g}_1\in G(\rho_0)= 1+m_{A_0}$. Comparing
top left entries we see that $\bar{g}_1=\bar{g}_2$ whence
$(g_1,g_2)\in G_1\times_{G_0}G_2=G_3$. This element conjugates
$\sigma$ to $\rho$ whence $[\rho]=[\sigma]$.
\end{pr}

This completes the proof of Thm. \ref{main}.

\section{Functoriality}

This section briefly illustrates that the usual functorial
properties of the universal deformation ring $R=R(G,\ol,N)$ hold
in our setting. Granting the universal property of $R$ this is
almost a formality.
\subsection{Morphisms}
Let $I$ be any set and consider the functor $A\mapsto
M_I(A)$ from $\C$ to topological $M_I(o)$-algebras. Suppose
$N\in\Modko$ with a pseudobasis indexed by $I$. Fix an isomorphism
$N\simeq k^{I}$. Given any morphism of functors
\[\delta: M_I\longrightarrow M_J\]
we may compose the diagram for $N$
\[\xymatrix{
    k[[H]] \ar[r]^-{cont}  & M_I(k)\\
  k[H] \ar[r]^{\subseteq} \ar[u]^{\subseteq} &  k[G] \ar[u].
  }
\] in the obvious way with $\delta(k)$ and obtain an element $N'\in\Modko$. Evidently $N'$ has
a pseudobasis indexed by $J$ which is why we refer to this
procedure as {\it change of range} for $N$. We denote by $\Irr$
the category consisting of the same objects as $\Modko$ but with
the following morphisms. Given $N,N'\in\Modko$ there is a morphism
$N\rightarrow N'$ if and only if $N'$ arises from $N$ via change
of range. In this situation a {\it morphism}
\[(G,\ol,N)\longrightarrow (G',\ol',N')\] consists by definition of
\begin{itemize}
\item[(a)] a morphism $\varphi: G'\rightarrow G$ of locally
$\Q_p$-analytic groups ({\it change of group}),

 \item[(b)] a local homomorphism
$\iota: \ol\rightarrow\ol'$ of commutative complete local
noetherian rings making $\ol'$ a pseudocompact $\ol$-algebra ({\it
change of base}),

\item[(c)] a morphism $N\rightarrow N'$ in $\Irr$ ({\it change of
range}).

\end{itemize}

The effect of such a morphism $(G,\ol,N)\rightarrow (G',\ol',N')$
on universal deformation rings is in complete analogy to the
compact case (cf. \cite{MazurI}, 1.3) which is why we omit the
details here.

\bigskip

Remark: To show the limits of this analogy let $k$ be finite and
consider the Pontryagin duality $\Modko\simeq\Modsmk$. The left
hand side inherits a duality coming from the smooth
contragredient. Due to the asymmetry of
$M_I(.)=\prod_I\oplus_I(.)$ with respect to "rows and columns" $N$
there does not seem to be a naive duality relating the universal
deformation rings of $N$ and its dual (comp. \cite{MazurI}, 1.3
(a.2)). For similar reasons there is no naive "determinant"
morphism (comp. [loc.cit.], 1.3 (a.3)).

\subsection{Tensor product}
We explain the effect of the tensor product on $\Modko$ on
universal deformation rings. Again, this is in analogy to loc.cit.
but since we will deduce the important proposition \ref{twisting}
from it we give some details. Let $N,N'\in\Modko$ and
$N'':=N\hato_k N'\in\Modko$. Given $A, A'\in\C$ with maximal
ideals $\m$ and $\m'$ respectively the results of sect.
\ref{pseudocompactrings} show that for any $A,A'\in\C$ the
$\ol$-algebra $A'':=A\hato_{\ol}A'$ lies again in $\C$ (and has
maximal ideal $\m''=\m\hato A'+A\hato\m'$). Given $[M,\al]$ and
$[M',\al]$ in $D_N(A)$ and $D_{\Nr'}(A')$ respectively the
compatibility of $\hato$ with direct products shows that
\[ M'':=M\hato_{\ol}M'\in {\rm Mod}^{\rm pro\,aug}_G(A'')^{\rm
fl}\]
and
\[M''/\overline{\m'' M''}\car M/\overline{\m M}\hato_{\ol}\,M'/\overline{\m' M'}\car N''\]
where the first isomorphism is canonical and the second induced by
$\al\otimes\al'$. Hence, $[M'',\al\otimes\al']\in D_{A''}(N'')$
whence a morphism of functors $D_N\times D_{N'}\rightarrow
D_{N''}$. If all functors involved are representable the usual
Yoneda lemma yields a homomorphism
\[h(\Nr,\Nr'): R''\longrightarrow R\hato_{\ol}R'\]
in $\C$ such that the system $(\Nr,\Nr')\mapsto h(\Nr,\Nr')$
satisfies the usual commutativity and associativity relations.
Now let $N,N'$ and $N''$ be as above but fix $[M,\al]\in
D_\Nr(\ol)$. We obtain a homomorphism ({\it contraction with a
lifting})
\[
R''\stackrel{h(\Nr,\Nr')}{\longrightarrow}R\hato_{\ol}R'\stackrel{M}{\longrightarrow}
R'\] satisfying the usual commutativity relations (loc.cit.).
It follows formally from these relations that if $|I|=1$ (so that
the $G$-action on $N$ is given by a character) this homomorphism
is always an isomorphism in $\C$ (the {\it twisting morphism} by
$M$) and satisfies the evident homomorphic property in the
variable $M$.

As we have observed before there is a canonical splitting $s:
k^\times\rightarrow\ol^\times$ in case $k$ is finite. Granting the
discussion above the following proposition follows therefore as in
loc.cit.
\begin{prop}\label{twisting}
Let $k$ be finite and suppose $D_N$ is representable. The
universal deformation ring $R(G,\ol,N)$ depends on $N$ only up to
twisting by characters.
\end{prop}

\section{Applications to $p$-adic Banach space
representations}\label{banach}

\subsection{Unitary Deformations}
We keep the above notations but assume that $\ol$ equals the ring
of integers in a finite extension $K$ of $\Q_p$. In particular,
$k$ is finite of characteristic $p>0$. Let $\varpi$ be a
uniformizer for $\ol$. Recall that a {\it Banach space
representation of $G$ over $K$} is a $K$-Banach space $V$ together
with a linear $G$-action such that the map
\[G\times V\longrightarrow V\]
giving the action is continuous (cf. \cite{STBanach}).
Together with continuous $G$-equivariant $K$-linear maps these
objects form a category ${\rm Ban}_G(K)$. We denote by ${\rm
Ban}_G(K)^{\leq 1}$ the category consisting of $K$-Banach spaces
$(V,||.||)$ such that $||V||\subseteq |K|$ endowed with a
continuous $G$-action such that $||.||$ is $G$-invariant (i.e.
$||gv||=||v||$ for all $g\in G,v\in V$). We let morphisms be
$G$-equivariant norm-decreasing $K$-linear maps. Elements
$(V,||.||)\in\Bangu$ will be called {\it unitary} representations.
\begin{lem}
Suppose $V$ is a $K$-Banach space representation of $G$ with a
$G$-invariant norm defining the topology. Then there is an
equivalent norm $||.||$ on $V$ which is $G$-invariant and such
that $||V||\subseteq |K|$.
\end{lem}
\begin{pr}
If $||.||'$ denotes the $G$-invariant norm on $V$ we put
\[||v||:=\inf \{r\in K: r\geq ||v||'\}\] for $v\in V$. Then $||.||$
is equivalent to $||.||'$ and $G$-invariant.
\end{pr}

Remark: It follows from this lemma that the full subcategory of
$\BanG$ consisting of elements $V$ that admit a $G$-invariant norm
is equivalent to the isogeny category of $\Bangu$.

\bigskip

After these preliminaries consider an element
$(V,||.||)\in\Bangu$. The unit ball $V^0:=\{v\in V:||v|\leq 1\}$
evidently inherits a $G$-action and
\[\bar{V}:=V^0/\varpi V^0\] defines a smooth $G$-representation
over $k$. We obtain in this way a functor
\[\Bangu\longrightarrow\Modsmk.\]
If $\bar{V}$ is admissible smooth $(V,||.||)$ is called {\it
admissible}.
\begin{theo}\label{defbanach}
Let $\rho$ be a smooth $G$-representation over $k$ which admits
only scalar endomorphisms (e.g admissible and absolutely
irreducible). Let $N=\rho^\vee$. There is a canonical and natural
bijection between the $\ol$-valued points of $R(G,\ol,N)$ and the
set of isomorphism classes of unitary $G$-Banach space
representations $V$ over $K$ such that $\bar{V}\simeq \rho$.
\end{theo}
\begin{pr}
Pick a compact open subgroup $H\subseteq G$. Given $V\in {\rm
Ban}_H(K)^{\leq 1}$ we may equip the $\ol$-module $V^d:={\rm
Hom}_\ol(V^0,\ol)$ with the topology of pointwise convergence and
the contragredient $H$-action. The $H$-equivariant version of the
discussion in \cite{STBanach}, (proof of) Thm. 1.2 shows that
$V\mapsto V^d$ establishes an equivalence of categories
\[(.)^d: {\rm Ban}_H(K)^{\leq 1}\car \PM (\ol[[H]])^{\rm
fl}.\]
It is evidently compatible with the restriction functors relative
to a compact open subgroup $H'\subseteq H$. Taking into account
the $G$-action it therefore restricts to an equivalence
\[(.)^d: {\rm Ban}_G(K)^{\leq 1}\car {\rm Mod}_G^{\rm pro\,aug}(o)^{\rm
fl}\] on the faithfully embedded subcategory ${\rm Ban}_G(K)^{\leq
1}$.
We may now form a diagram of functors


\begin{equation}\label{diag}\xymatrix{
{\rm Ban}_G(K)^{\leq 1} \ar[r]^{\sim}_{(.)^d}\ar[d]^{V\mapsto
\bar{V}} & {\rm Mod}_G^{\rm pro\,aug}(o)^{\rm
fl} \ar[d]^{\pi^*}\\
  {\rm Mod}^{\rm sm}_G(k)  \ar[r]^{\sim}_{(.)^\vee}& {\rm Mod}_G^{\rm pro\,aug}(k)}
\end{equation}

where the right perpendicular arrow refers to the base change
relative to the residue map $\pi: \ol\rightarrow k$ and the lower
horizontal arrow equals Pontryagin duality. A straightforward
equivariant version of \cite{Paskunas}, Lem. 5.4 proves the
diagram to be commutative. It follows that the functor $(.)^d$
induces a natural and canonical bijection between the set of
isomorphism classes in $\Bangu$ of $(V,||.||)$ with $\bar{V}\simeq
\rho$ and $D_\Nr(\ol)$.
\end{pr}

\subsection{Principal series representations of
${\rm GL}_2(\Q_p)$}

To illustrate the above methods we let $G={\rm
GL}_2(\mathbb{Q}_p)$ and compute the fibers of the reduction map
in principal series representations. Besides deformation theory
our computations strongly rely on results of M. Emerton concerning
the functor of {\it ordinary parts} (cf.
\cite{EmertonI},\cite{EmertonII}). To start with let $P$ and
$\overline{P}$ be the Borel subgroup of $G$ consisting of upper
triangular and lower triangular matrices respectively. Given two
smooth characters $\chi_i: \Q_p^\times\longrightarrow A^\times,
i=1,2$ for $A\in C$ we view $\chi=\chi_1\otimes\chi_2$ as a smooth
character of the diagonal torus $T$ in $G$ in the obvious way.
Define
\[ {\rm Ind}_{\oP}^G(\chi)=\{f: G\rightarrow A | f {\rm ~locally~constant},
f(\overline{p}g)=\chi(\overline{p})f(g), p\in\oP,g\in G\}\] with
$G$-action by right translations. It is a smooth admissible
$G$-representation over $A$ in the sense of \cite{EmertonI}, Def.
2.2.5. Finally, let $\epsilon(a)=a|a|\in\Z_p^\times$ for all
$a\in\Q_p^\times$ and write $\bar{\epsilon}$ for the induced
smooth character $\Q_p^\times\rightarrow k^\times$.

If $V^0$ denotes the unit ball of an element $(V,||.||)$ in ${\rm
Ban}_G(K)^{\leq 1}$ we write $V_n:=V^0/\varpi^nV^0$ and
$\ol_n:=\ol/\varpi^n\ol$ for all $n$.
\begin{lem}\label{faith} The $\ol_n$-module $V_n$ is faithful for all $n$.
\end{lem}
\begin{pr}
By the same argument as in case $n=1$ the diagram (\ref{diag})
remains commutative when we replace $k$ by $\ol_n$, $\varpi$ by
$\varpi^n$ and restrict to topologically free $\ol_n$-modules in
the lower right corner. But then $(V_n)^\vee$ is topologically
free and therefore $V_n$ is faithful.
\end{pr}

After these preliminaries let us fix smooth characters $\ochi_i:
\Q_p^\times\longrightarrow k^\times, i=1,2$ for which
$\bar{\chi}_1\bar{\chi}_2^{-1}\neq 1,\oep$. The $G$-representation
${\rm Ind}_{\oP}^G(\bar{\chi})$ is then admissible and absolutely
irreducible (cf. \cite{BarthelLivne}). In particular, it admits
only scalar endomorphisms (cf. Lem. \ref{schur1}). Suppose our
chosen Banach space representation $V$ satisfies
\[\oV=V^0/\varpi V^0\simeq {\rm Ind}_{\oP}^G(\bar{\chi}).\]
\begin{lem}\label{deformlemma}
Let $n\geq 1$ and suppose there is an isomorphism $V_n\car {\rm
Ind}_{\oP}^G(\ochi^n)$ with some smooth $\ol_n^\times$-valued
character $\ochi^{(n)}=\chi^{(n)}_1\otimes\chi^{(n)}_2$. Then
there is a smooth $\ol_{n+1}^\times$-valued character
$\ochi^{(n+1)}=\chi^{(n+1)}_1\otimes\chi^{(n+1)}_2$ and a
commutative diagram of smooth $G$-representations
\[\xymatrix{
  V_{n+1} \ar[r]^-{\sim}_-{\varphi_{n+1}}\ar[d]^-{mod\;\varpi^n} &  {\rm
Ind}_{\oP}^G(\ochi^{n+1}) \ar[d]^-{mod\;\varpi^{n}}\\
  V_n \ar[r]^-{\sim}_-{\varphi_n}&  {\rm
Ind}_{\oP}^G(\ochi^{n}).   }
\]
\end{lem}
\begin{pr}
Granting the above lemma this is a straightforward generalization
of \cite{EmertonII}, Prop. 4.1.5. 
\end{pr}

On the other hand, given a continuous character $\chi:
T\rightarrow \ol^\times$, we may define
\[ {\rm ^cInd}_{\oP}^G(\chi)=\{f: G\rightarrow K | f
{\rm ~continuous}, f(\overline{p}g)=\chi(\overline{p})f(g),
p\in\oP,g\in G\}\] with $G$ acting by right translations. Equipped
with the supremum norm taken over the compact space $\oP\backslash
G$ it constitutes an admissible unitary $G$-Banach space
representation over $K$, the so-called {\it ordinary continuous
principal series}. Its irreducibility properties are well-known
(cf. \cite{EmertonL}, Prop. 5.3.4). The reduction mod $\varpi$ of
its unit ball is of the form ${\rm
Ind}_{\bar{P}}^{G}(\bar{\chi})$. Recall that $\kappa=2$ if $p=2$
and $1$ else. Let $\mu'$ denote the group with $\kappa$ elements.
\begin{theo}
If $V\in {\rm Ban}_G(K)^{\leq 1}$ such that $\oV\simeq {\rm
Ind}_{\oP}^G(\bar{\chi})$ then $V={\rm ^cInd}_{\oP}^G(\chi)$ with
a lift $\chi$ of $\ochi$ to $\ol^\times$. The isomorphism classes
of such $V$ are therefore in bijection to the product of four
copies of the maximal ideal $(\varpi)$ with two copies of $\mu'$.
\end{theo}
\begin{pr}
The above lemma together with \cite{EmertonI}, Lem. 4.1.1 implies
that the injective map \[\chi\mapsto {\rm ^cInd}_{\oP}^G(\chi)\]
from unitary lifts of $\bar{\chi}$ to unitary lifts of $\oV$ is
surjective. Using the above theorem it remains to compute
$D_N(\ol)$ where $N=\bar{\chi}^\vee=\bar{\chi}^{-1}$. But the
example after Cor. \ref{noetherian} shows that
\[R(T,\ol,N)=\ol[[x_1,x_2,x_3,x_4]]\otimes_\ol\ol[\mu'\times\mu'].\]
\end{pr}

Remark: Let $H^\bullet{\rm Ord}_P$ denote the $\delta$-functor
associated to the functor of ordinary parts ${\rm Ord}_P$ relative
to $P$. The proof of Lem. \ref{deformlemma} relies on the fact
that \[H^1{\rm Ord}_P ({\rm Ind}_{\oP}^G(\bar{\chi}))\simeq
\bar{\chi}_2\bar{\epsilon}^{-1}\otimes\bar{\chi}_1\bar{\epsilon}\]
as smooth $T$-representations (cf. \cite{EmertonII}, Thm. 4.1.3
(1)). A generalization of the above corollary to other groups than
${\rm GL}_2(\Q_p)$ is therefore related to a better understanding
of (higher) ordinary parts.

\bibliography{deformbib}

\providecommand{\bysame}{\leavevmode\hbox to3em{\hrulefill}\thinspace}
\providecommand{\MR}{\relax\ifhmode\unskip\space\fi MR }
\providecommand{\MRhref}[2]{%
  \href{http://www.ams.org/mathscinet-getitem?mr=#1}{#2}
}
\providecommand{\href}[2]{#2}
\begin{thebibliography}{DdSMS99}

\bibitem[BL94]{BarthelLivne}
L.~Barthel and R.~Livn{\'e}, \emph{Irreducible modular representations of
  {${\rm GL}\sb 2$} of a local field}, Duke Math. J. \textbf{75} (1994), no.~2,
  261--292.

\bibitem[B{\"o}c98]{Boeckle}
G.~B{\"o}ckle, \emph{The generic fiber of the universal deformation space
  associated to a tame galois representation}, Manuscripta Math. \textbf{96}
  (1998), 231--246.

\bibitem[Bou89]{B-GT}
Nicolas Bourbaki, \emph{General topology. {C}hapters 5--10}, Elements of
  Mathematics (Berlin), Springer-Verlag, Berlin, 1989, Translated from the
  French, Reprint of the 1966 edition.

\bibitem[Bre]{BreuilNYC}
C.~Breuil, \emph{Representations of {G}alois and of {${\rm GL}_2$} in
  characteristic {$p$}}, Graduate course at Columbia University (fall 2007).
  Available at: http://www.ihes.fr/{$\sim$}breuil/PUBLICATIONS/New-York.pdf.

\bibitem[Bru66]{Brumer}
Armand Brumer, \emph{Pseudocompact algebras, profinite groups and class
  formations}, J. Algebra \textbf{4} (1966), 442--470.

\bibitem[Col10]{Colmez10}
P.~Colmez, \emph{Repr\'esentations de {${\rm GL}\sb 2(\bold Q\sb p)$} et
  {$(\phi,\Gamma)$}-modules}, Ast\'erisque (2010), no.~330, 281--509.
  \MR{2642409}

\bibitem[DdSMS99]{DDMS}
J.~D. Dixon, M.~P.~F. du~Sautoy, A.~Mann, and D.~Segal, \emph{Analytic
  pro-{$p$} groups}, second ed., Cambridge Studies in Advanced Mathematics,
  vol.~61, Cambridge University Press, Cambridge, 1999.

\bibitem[Eme]{EmertonA}
M.~Emerton, \emph{Locally analytic vectors in representations of locally
  $p$-adic analytic groups}, Preprint. To appear in: Memoirs of the AMS.

\bibitem[Eme06]{EmertonL}
\bysame, \emph{A local-global compatibility conjecture in the {$p$}-adic
  {L}anglands programme for {${\rm GL}\sb {2/{\Bbb Q}}$}}, Pure Appl. Math. Q.
  \textbf{2} (2006), no.~2, part 2, 279--393.

\bibitem[Eme10a]{EmertonI}
\bysame, \emph{Ordinary parts of admissible representations of $p$-adic
  reductive groups {I}. {D}efinitions and first properties}, Ast\'erisque
  \textbf{331} (2010), 335--381.

\bibitem[Eme10b]{EmertonII}
\bysame, \emph{Ordinary parts of admissible representations of $p$-adic
  reductive groups {II}.{D}erived functors}, Ast\'erisque \textbf{331} (2010).

\bibitem[Gab62]{GabrielPhD}
P.~Gabriel, \emph{Des {C}at\'egories {A}beli\'ennes}, Bull. {S}oc. {M}ath.
  {F}rance \textbf{90} (1962), 323--448.

\bibitem[Gab70]{SGA3}
\bysame, \emph{{\'E}tude infinit\'esimale des sch\'emas en groupes. {I}n:
  {S}ch\'emas en groupes {I}, {E}xp. {VIIB}}, Lect. Notes Math., vol. 151,
  Springer, Berlin-Heidelberg-New York, 1970.

\bibitem[Gou01]{GouDef}
Fernando~Q. Gouv{\^e}a, \emph{Deformations of {G}alois representations},
  Arithmetic algebraic geometry ({P}ark {C}ity, {UT}, 1999), IAS/Park City
  Math. Ser., vol.~9, Amer. Math. Soc., Providence, RI, 2001, Appendix 1 by
  Mark Dickinson, Appendix 2 by Tom Weston and Appendix 3 by Matthew Emerton,
  pp.~233--406.

\bibitem[Gro]{GrothendieckDesII}
A.~Grothendieck, \emph{Technique de descente et th\'eor\`emes d'existence en
  g\'eom\'etrie alg\'ebrique. {II}. {L}e th\'eor\`eme d'existence en th\'eorie
  formelle des modules}, S\'eminaire {B}ourbaki, {V}ol.\ 5, Soc. Math. France,
  Paris, pp.~Exp.\ No.\ 195, 369--390.

\bibitem[Kis10]{KisinDEF}
M.~Kisin, \emph{{D}eformations of {$G_{\bf {Q}_p}$} and {${\rm
  GL}_2(\mathbb{Q}_p)$} representations.}, Ast\'erisque \textbf{330} (2010),
  529--542.

\bibitem[Maz89]{MazurI}
Barry Mazur, \emph{Deforming {G}alois representations}, Galois groups over
  {${\bf Q}$} ({B}erkeley, {CA}, 1987), Math. Sci. Res. Inst. Publ., vol.~16,
  Springer, New York, 1989, pp.~385--437.

\bibitem[Maz97]{MazurII}
\bysame, \emph{An introduction to the deformation theory of {G}alois
  representations}, Modular forms and {F}ermat's last theorem ({B}oston, {MA},
  1995), Springer, New York, 1997, pp.~243--311.

\bibitem[Nag62]{Nagata}
Masayoshi Nagata, \emph{Local rings}, Interscience Tracts in Pure and Applied
  Mathematics, No. 13, Interscience Publishers a division of John Wiley \&
  Sons\, New York-London, 1962.

\bibitem[Pas10]{Paskunas}
V.~Paskunas, \emph{Admissible unitary completions of locally {$\Q_p$}-rational
  representations of {${\rm GL}_2(F)$}}, Representation {T}heory \textbf{14}
  (2010), 324--354.

\bibitem[RZ00]{Ribes}
Luis Ribes and Pavel Zalesskii, \emph{Profinite groups}, Ergebnisse der
  Mathematik und ihrer Grenzgebiete. 3. Folge. A Series of Modern Surveys in
  Mathematics, vol.~40, Springer-Verlag, Berlin, 2000.

\bibitem[Sch68]{Schlessinger}
Michael Schlessinger, \emph{Functors of {A}rtin rings}, Trans. Amer. Math. Soc.
  \textbf{130} (1968), 208--222.

\bibitem[ST02]{STBanach}
P.~Schneider and J.~Teitelbaum, \emph{Banach space representations and
  {I}wasawa theory}, Israel J. Math. \textbf{127} (2002), 359--380.

\bibitem[SV06]{SV1}
Peter Schneider and Otmar Venjakob, \emph{On the codimension of modules over
  skew power series rings with applications to {I}wasawa algebras}, J. Pure
  Appl. Algebra \textbf{204} (2006), no.~2, 349--367.

\bibitem[vGvdB97]{GvB}
M.~van Gastel and M.~van~den Bergh, \emph{Graded modules of gelfand-kirillov
  dimension one over three-dimensional artin-schelter regular algebras}, J.
  Algebra \textbf{196} (1997), 251--282.

\bibitem[Wil98]{Wilson}
John~S. Wilson, \emph{Profinite groups}, London Mathematical Society
  Monographs. New Series, vol.~19, The Clarendon Press Oxford University Press,
  New York, 1998.

\end{thebibliography}
\bibliographystyle{amsalpha}

\end{document}